\documentclass[11pt,twoside,a4paper,leqno]{article}

\usepackage[left=2.5cm,right=2.5cm,top=3cm,height=24cm,headheight=4cm, headsep=4mm, footskip =1cm]{geometry}

\usepackage{mathrsfs,amsthm,amsfonts,amssymb}

\newcommand {\R}{\mathbb{R}} 
\newcommand {\C}{\mathbb{C}} 

\newtheorem{theo}{Theorem}
\newtheorem{lem}{Lemma}

\usepackage{fancyhdr}
\begin{document}
\thispagestyle{empty}
\begin{center}
\begin{LARGE}
\textbf{A gradient estimate for the Hermitian}

\vskip1mm

\textbf{Monge-Ampère equation}

\vskip1.2cm

\textbf{Abdellah HANANI}{\footnote {\noindent Current address : Université de Lille, département mathématiques, Bât. M2, 59655, Villeneuve d’Ascq
Cedex, France.\\ email address: abdellah.hanani@univ-lille.fr}}
\end{LARGE}
\end{center}

\vskip6mm

\noindent\rule{\linewidth}{1pt}

\vskip4mm

\noindent\textbf{Abstract.} We improve our previous gradient estimate for the Monge-Ampère equation on a compact Hermitian manifold and give a estimates for the non-mixed second order derivatives. These estimates are required to apply either the Evans-Krylov $C^{2,\alpha}$ estimates or the third derivatives estimates for equations with a gradient term.

\vskip6mm

\noindent\textbf{Keyword.} Monge-Amp\`ere equation, Hermitian manifold, gradient estimate.
      
\vskip2mm

\noindent Mathematical Subject Classification : 32U05, 32U40, 53C55

\vskip2mm

\noindent\rule{\linewidth}{1pt}

\vskip1cm

\section{Introduction}

\vskip4mm

\noindent The complex Monge-Amp\`ere equation over the last half century has been the subject of intensive investigations. In the framework of K\"ahler geometry, a source of interest in this equation goes back to the well known conjecture of Calabi : To every closed (1,1)-form $R$ representing the first Chern class of a Kähler manifold $(M,g)$ of complex dimension $n\geq 2$, there is a unique Kähler metric in the same Kähler class whose Ricci form is $R$.

\vskip2mm

\noindent In case the complex manifold has vanishing first Chern class, write
$$\omega =\mathrm{i}g_{\lambda \overline{\mu}}\mathrm{d}z^{\lambda}\wedge \mathrm{d}z^{\overline{\mu}},$$
Calabi realized that his conjecture can be reduced to solving a Monge-Ampère equation of the following form
$$(\omega+\mathrm{i}\partial\overline{\partial}\varphi)^n=e ^f\omega ^n\quad \mbox{with}\quad \omega+\mathrm{i}\partial\overline{\partial}\varphi>0,$$
for a real valued function $\varphi$, where $f$ is a smooth real function on $M$. The equation was solved by Aubin [1-2] and Yau [13]. The main purpose of this paper is to go back on those estimates which were not required in the above authors works but became of considerable interest in the last two decades.

\vskip3mm

\noindent The two authors were able to derive independently estimates for the Laplacian of the solution assuming only a uniform estimate of the solution. The gradient estimate could be bypassed. On the other hand, the third derivatives estimate, given once more independently by the two authors, did not make use of the estimate of second order derivatives of pure indexes. To deal with equations involving a gradient term, we need to derive these bypassed estimates.

\vskip3mm

\noindent Let $(M,g)$ be a compact Hermitian manifold, without boundary. Denote by $\nabla$ the Chern connection and by $TM$ the tangent bundle of $M$. For a nonnegative smooth real valued function $F$ on $TM\times \R$, we consider the following Monge-Amp\`ere equation
\begin{equation}
\left(\omega+\mathrm{i}\partial\overline{\partial}\varphi\right)^n=F\omega ^n,
\end{equation}
with the condition that $\omega+\mathrm{i}\partial \overline{\partial}\varphi>0$. The function $\varphi$ is said to be admissible. For $(x,p)\in TM$ and $t\in \R$, we write 
$$F_x(\bullet )=F(\bullet ,p,t),\quad F_p(\bullet )=F(x,\bullet ,t)\quad \mbox{and}\quad F_t(\bullet )=F(x,p,\bullet ).$$
Our gradient estimate may be stated as follows.
 
\vskip6mm

\begin{theo} Let $\varphi \in C^{3}(M)$ be an admissible solution of $(1)$ such that $\Vert \varphi \Vert _{\infty}\leq C_0$ for a positive constant $C_0$. Suppose that there exists a real number $r$ such that $r<2$ and that, for all $(x,p,t)\in TM\times [-C_0,C_0]$,
\begin{equation}
F\leq C\left(1+\vert p\vert ^r\right) \quad \mbox{and}\quad \left|\nabla \left[(F^{\frac{1}{n}})_x\right]\right|,\; \left|\left[(F^{\frac{1}{n}})_t\right]'\right|,\; \left|p^{\alpha}\frac{\partial (F^{\frac{1}{n}})_p}{\partial p^{\alpha}}\right| \leq CF^{\frac{1}{n}},
\end{equation}
where $(p^a)$ are the natural coordinates in the fibers of $TM$. Then 
$$\vert \nabla \varphi \vert \leq C_1,$$
where the constant $C_1$ depends only on $C_0$, $M$ and the constant in $(2)$. In case $F$ does not depend on the gradient, hypothesis (2) is not required and $C_1$ depends only on $C_0$, $M$ and an upper bound for $\Vert \nabla F^{\frac{1}{n}}\Vert _{\infty}$.
\end{theo}

\vskip2mm

\noindent This gradient estimate appears for the first time in [8] and [9] in the nineties. For the next few years there was no activity on the Monge-Ampère equation on Hermitian manifolds until very recently, when the results were rediscovered P. Guan [7], B. Guan-Q. Li [6], D.H Phong-J. Sturm [11] and . Moreover, S.Dinew-S.Kolodziej [4] also studied the equation in the weak sense and obtained the uniform estimates of the solution via suitably constructed pluripotential theory.

\vskip4mm

\noindent Following the same strategy as that of Aubin and Yau, see also P. Cherrier [3], one can now estimate the Laplacian. It follows from this estimate that the metrics $g$ and $g'$ defined by $(g_{\lambda \overline{u}}+\partial _{\lambda \overline{u}}\varphi)$ are equivalent : there exists a strictly positive constant $a$ such that
$$a^{-1}g\leq g'\leq ag.$$
The dependency of the constant $a$ on the geometry of $M$ is quite explicit. It also depends on a $C^2$-norm of $F$ on the compact
$${\mathscr K}=\{(x,p,t)\in TM\times \R\;:\; x\in M,\ \vert p\vert \leq C_1,\
\vert t\vert \leq C_0\}.$$

\noindent The next step in Aubin and Yau works was to estimate the third derivatives. For equations with gradient term, one can not handle this estimate directly without completing the $C^2$-estimate. These estimates are also required to apply the Evans-Krylov theorem for $C^{2,\alpha}$-estimates. We have the following estimates for the second derivatives of non-mixed type.

\vskip6mm

\begin{theo} Let $\varphi \in \C^{4}(M)$ be an admissible
solution of equation $(1)$ such that $\|\varphi\|\leq C_0$ and $\|\nabla \varphi\|\leq C_1$. Then
$$\vert \nabla _{\lambda \mu}\varphi\nabla ^{\lambda \mu}\varphi\vert \leq C_2,$$
where the constant $C_2$ depends only on $C_0$, $C_1$, on $M$ in terms of its geometry, a bound of the $C^{2}$-norm of $F$ on the compact ${\mathscr K}$ and $\displaystyle \min _{\mathscr K}F$.
\end{theo}

\vskip8mm

\section{Proof of theorem 1}

\vskip4mm

In local coordinates $(z^{\lambda})$, write  $\displaystyle z^{\overline{\lambda}}=\overline{z^{\lambda}}$, $\displaystyle e_{\lambda}=\partial _{\lambda}=\frac{\partial}{\partial z^{\lambda}}$,
$$g=g_{\lambda \overline{\mu}}\mathrm{d} z^{\lambda}\otimes \mathrm{d}z^{\overline{\mu}}\quad \mbox{and}\quad g^{-1}=g^{\lambda \overline{\mu}}e _{\lambda}\otimes e_ {\overline{\mu}}.$$
The Laplacians $\Delta$ and $\Delta'$ with respect to the metrics $g$ and $g'$ acting on functions are defined as follows
$$\Delta \varphi =g^{\alpha \overline{\beta}}\nabla _{\alpha \overline{\beta}}\varphi \quad \mbox{and}\quad \Delta '\varphi =g'^{\alpha \overline{\beta}}\nabla _{\alpha \overline{\beta}}\varphi.$$
The norms of a tensor $T=T_{\alpha}\mathrm{d} z^{\alpha}$ with respect to the metrics $g$ and $g'$ are denoted as follow $g$ and $g'$ acting on functions are defined as follows
$$|T|^2=g^{\alpha \overline{\beta}}T_{\alpha}T_{ \overline{\beta}} \quad \mbox{and}\quad |T|'^2=g'^{\alpha \overline{\beta}}T_{\alpha}T_{\overline{\beta}} .$$
Since, we are restricting ourselves to the smooth case, we can write equation (1) as follows
\begin{equation}
\det \left(\delta ^{\lambda}_{\mu}+\nabla ^{\lambda}_{\mu}\varphi \right)=F(x,\nabla \varphi ,\varphi).
\end{equation}
A key ingredient in the proof of theorem 1 is the following lemma.

\vskip6mm

\begin{lem} Let $\varphi \in C^{3}(M)$ be an admissible
solution of $(3)$. Then
\begin{equation}
\begin{array}{ccc}\Delta '\vert \nabla \varphi \vert ^{2} &=& 2\mathrm{Re}\left<\nabla \log F,\nabla \varphi \right>_g+g'^{\alpha \overline \beta }\nabla _{\alpha \lambda  }\varphi \nabla _{\overline \beta }^{\ \ \lambda }\varphi +g'^{\alpha \overline\beta }g_{\alpha \overline\beta}-n+\Delta \varphi \\[3mm]&&+2\mathrm{Re}\left(g'^{\alpha \overline{\beta }}T^{\nu }_{\lambda \alpha }\nabla ^{\lambda }\varphi \nabla _{\nu \overline{\beta}}\varphi \right)+g'^{\alpha \overline{\beta }}R_{\alpha\overline{\beta }\rho\overline{\mu}}\nabla ^{\rho}\varphi \nabla ^{\overline{\mu }}\varphi ,\end{array}
\end{equation}
The notation $\mathrm{Re}(z)$ stands for the real part of the complex $z$.
\end{lem}

\vskip4mm

\noindent\textit{Proof.} In the following, we denote $A=g'^{\alpha \overline \beta }\nabla _{\alpha \lambda  }\varphi \nabla _{\overline \beta }^{\lambda }\varphi$. Let us compute $\Delta '\varphi $,
$$\Delta '\vert  \nabla \varphi \vert ^{2}=2\mathrm{Re}(g'^{\alpha \overline \beta}\nabla _{\alpha \overline \beta \lambda }\varphi \nabla ^{\lambda }\varphi )+g'^{\alpha \overline \beta }\nabla _{\alpha \lambda  }\varphi \nabla _{\overline \beta }^{\lambda }\varphi+g'^{\alpha\overline \beta }\nabla _{\alpha }^{\lambda}\varphi \nabla _{\lambda \overline \beta }\varphi.$$
To interchange the order of derivation in the third derivative terms, we use the following relations
\begin{equation}
\nabla _{\alpha \overline\beta \lambda }\varphi =\nabla _{\lambda \alpha  \overline{\beta }}\varphi +T^{\nu }_{\lambda \alpha }\nabla _{\nu \overline\beta }\varphi 
\end{equation}
and

\begin{equation}
\nabla _{\alpha \overline{\beta \mu }}\varphi =\nabla _{\overline{\mu }\alpha \overline{\beta }}\varphi +T^{\overline\nu}_{\overline{\mu \beta }}\nabla _{\alpha \overline{\nu }}\varphi +g^{\rho \overline \nu
}R_{\alpha \overline{\beta }\rho \overline{\mu }}\nabla _{\overline{\nu }}\varphi .
\end{equation}
So that
\begin{equation}
\Delta '\vert  \nabla \varphi \vert ^{2}=2\mathrm{Re}(g'^{\alpha \overline \beta}\nabla _{\lambda\alpha \overline \beta }\varphi \nabla ^{\lambda }\varphi )+g'^{\alpha \overline \beta }\nabla _{\alpha \lambda  }\varphi \nabla _{\overline \beta }^{\lambda }\varphi+g'^{\alpha\overline \beta }\nabla _{\alpha }^{\lambda}\varphi \nabla _{\lambda \overline \beta }\varphi +E.
\end{equation}
where
$$\displaystyle E=2\mathrm{Re}(g'^{\alpha \overline{\beta }}T^{\nu }_{\lambda \alpha }\nabla ^{\lambda }\varphi \nabla _{\nu \overline{\beta}}\varphi )+g'^{\alpha \overline{\beta }}g^{\rho \overline{\nu}}R_{\alpha\overline{\beta }\rho\overline{\mu}}\nabla
_{\overline\nu }\varphi \nabla ^{\overline{\mu }}\varphi .$$
Using the definition of $g'$, we see that 
\begin{equation}
\nabla _{\nu \overline\beta}\varphi=g'_{\nu  \overline\beta}-g_{\nu \overline\beta}\quad \mbox{ and }\quad g'^{\alpha \overline\beta}g'_{\nu  \overline\beta}=\delta
^{\alpha}_{\nu}.
\end{equation}
So
\begin{equation}
g'^{\alpha\overline \beta }\nabla _{\alpha }^{\lambda}\varphi \nabla _{\lambda \overline \beta }\varphi=+g'^{\alpha \overline\beta }g_{\alpha \overline\beta}-n+\Delta \varphi .
\end{equation}
Now, differentiating once equation $(3)$, we obtain
\begin{equation}
g'^{\alpha \overline{\beta  }}\nabla _{a\alpha  \overline{\beta }}\varphi =\nabla  _{a}\log f \quad \mbox{for all }a\in \{1,\dots ,n,\overline{1},\dots ,\overline{n}\}.
\end{equation}
Reporting (9) and (10) into (7) gives the relation (4).

\vskip4mm

\noindent\textbf{Proof of theorem 1.} Let $k$ be strictly positive real to be determined later. Denote
$$\Gamma (\varphi)=\vert\nabla \varphi \vert  ^{2}\exp[\exp[k(C_0-\varphi )]].$$
Since $M$ is compact, $\Gamma $ attains its maximum at some point $P\in M$. We are done if we prove that $\vert \nabla \varphi \vert (P)\leq C$. So, assume that
\begin{equation}
\vert \nabla \varphi \vert (P)\geq 1.
\end{equation}
The function $\log \Gamma $ attains also a maximum at $P$. Therefore
\begin{equation}
\frac{\nabla _{\alpha }\vert  \nabla  \varphi \vert ^2}{\vert \nabla \varphi \vert  ^{2} }-ke^{k(C_0-\varphi )}\nabla  _{\alpha }\varphi =0,
\end{equation}
and writing that the Laplacian $\Delta '\log \Gamma =g'^{\alpha \overline \beta}\nabla _{\alpha \overline \beta}\log \Gamma $ of $\log \Gamma $ is non positive at $P$, we obtain
\begin{equation}
\frac{\Delta '\vert \nabla  \varphi \vert ^2}{| \nabla  \varphi |^2}-\frac{{|\nabla |\nabla \varphi |^2|'}^2}{|\nabla  \varphi |^4}+ke^{k(C_0-\varphi )}(k\vert \nabla \varphi \vert '^{2}-\Delta '\varphi )\leq 0.
\end{equation}

\vskip2mm

\noindent We treat first the second term in the left hand side of $(13)$. Since $\displaystyle {\vert \nabla \log\Gamma (\varphi )\vert }'^{2}(P)=0$, we can write
\begin{equation}
\frac{{\vert \nabla \vert \nabla  \varphi \vert  ^2\vert '}^2 }{\vert \nabla \varphi \vert ^4}=-k^{2}e^{2k(C_0-\varphi )}{\vert \nabla \varphi \vert '}^{2}+2\frac{ke^{k(C_0-\varphi )} }{\vert \nabla \varphi \vert ^2}\mathrm{Re}(g'^{\alpha \overline\beta }\nabla _{\alpha }\vert \nabla \varphi \vert ^{2}\nabla _{\overline{\beta }}\varphi ).
\end{equation}
Expanding $\displaystyle \nabla \vert\nabla\varphi \vert 
^{2}$, in view of $(8)$, we see that
$$2\mathrm{Re}\left (g'^{\alpha \overline{\beta }}\nabla _{\alpha }\vert \nabla \varphi \vert  ^{2}\nabla  _{\overline{\beta }}\varphi\right)=2\left(\vert\nabla \varphi \vert  ^{2}-\vert  \nabla \varphi \vert '^{2}\right)+2\mathrm{Re}\left(g'^{\alpha  \overline{\beta }}g^{\lambda \overline{\mu  }}\nabla _{\alpha  \lambda  }\varphi \nabla _{\overline{\mu }}\varphi \nabla
 _{\overline{\beta }}\varphi \right).$$
Therefore, using the inequality of Cauchy, we get 
$$2\mathrm{Re}\left(g'^{\alpha \overline{\beta  }}g^{\lambda  \overline{\mu}}\nabla _{\alpha \lambda  }\varphi \nabla  _{\overline{\mu }}\varphi\nabla _{\overline{\beta }}\varphi\right)\leq \frac{A}{ke^{k(C_0-\varphi )} }+ke^{k(C_0-\varphi)}\vert \nabla \varphi\vert  ^{2}\vert \nabla \varphi \vert '^{2}$$
and reporting into (14) gives us the desired control. Namely,
\begin{equation}
\frac{\vert \nabla '\vert \nabla \varphi \vert ^{2}\vert ^{2} }{\vert \nabla \varphi \vert ^4}\leq \frac{A}{\vert \nabla \varphi \vert ^2}+2ke^{k(C_0-\varphi
)}.
\end{equation}

\vskip2mm

\noindent To deal with the first term in right hand side of $(13)$, notice that
$$\frac{1}{n}\nabla _{\lambda}\log F=\frac{\nabla _{\lambda}F^\frac{1}{n}}{F^{\frac{1}{n}}}=\frac{\nabla _{\lambda}\left[(F^{\frac{1}{n}})_x\right]}{F^{\frac{1}{n}}}+\frac{\left[(F^{\frac{1}{n}})_t\right]'}{F^{\frac{1}{n}}}\nabla _{\lambda}\varphi +\frac{\frac{\partial (F^{\frac{1}{n}})_p}{\partial p^{a}}}{F^{\frac{1}{n}}}\nabla _{\lambda a}\varphi .$$
So that, taking account of the relations
$$\nabla _{\lambda \overline \mu}\varphi=\nabla _{\overline \mu\lambda }\varphi\quad \mbox{and}\quad \nabla _{\lambda \mu}\varphi=\nabla _{\mu
\lambda }\varphi+T^{\nu}_{\mu \lambda }\nabla _{\nu}\varphi,$$
we obtain the following equality
$$\begin{array}{ccl}\displaystyle \frac{2}{n}\mathrm{Re}\left(\nabla ^{\lambda  }\varphi \nabla _{\lambda }\log F\right)&=&\displaystyle \frac{2}{F^{\frac{1}{n}}}\mathrm{Re}\left(\sum _{\alpha =1}^n\frac{\partial (F^{\frac{1}{n}})_p}{\partial p^{\alpha }}\nabla _{\alpha}\vert \nabla \varphi \vert ^{2}\right)+\frac{2}{F^{\frac{1}{n}}}\mathrm{Re}\left(\nabla ^{\lambda }\varphi \nabla _{\lambda }(F^{\frac{1}{n}})_x\right)\\[4mm] & &\displaystyle +\frac{2(F^{\frac{1}{n}})'_t\vert \nabla \varphi \vert }{F^{\frac{1}{n}}}+f'_t\vert \nabla \varphi \vert \mathrm{Re}\left(\sum _{\lambda =1}^n\frac{\partial (F^{\frac{1}{n}})_p}{\partial p^{\lambda }}T^{\nu}_{\lambda \mu}\nabla ^{\mu}\varphi \nabla _{\nu}\varphi \right),\end{array}$$
Now, in view of (12), we can write
$$\mathrm{Re}\left(\sum _{\alpha =1}^n\frac{\partial (F^{\frac{1}{n}})_p}{\partial p^{\alpha }}\nabla _{\alpha}\vert \nabla \varphi \vert ^{2}\right)=k e^{k(C_0-\varphi)}|\nabla \varphi|^2\mathrm{Re}\left(\sum _{\alpha =1}^n\frac{\partial f_p}{\partial p^{\alpha }}\nabla _{\alpha}\varphi\right)$$
and the assumption (11) and the structure conditions $(2)$ imply that
\begin{equation}
2\mathrm{Re}(\nabla  ^{\lambda  }\varphi \nabla _{\lambda  }\log F)\geq -C_3k e^{k(C_0-\varphi)}|\nabla \varphi|^2,
\end{equation}
where the constant $C_3$ depends on $C_0$, the constant in (2) and $M$.

\vskip2mm

\noindent Inserting inequality $(16)$ into $(4)$ and taking account of the fact that
$$2\mathrm{Re}\left(g'^{\alpha \overline{\beta }}T^{\nu }_{\lambda \alpha }\nabla ^{\lambda }\varphi \nabla _{\nu \overline{\beta}}\varphi \right)+g'^{\alpha \overline{\beta }}R_{\alpha\overline{\beta }\rho\overline{\mu}}\nabla ^{\rho}\varphi \nabla ^{\overline{\mu }}\varphi \geq -C_4\vert \nabla \varphi \vert  ^{2}(1+g'^{\alpha\overline\beta }g_{\alpha \overline\beta }),$$
where the constant $C_4$ depends only on an upper bound on $M$ of the norm of the torsion tensor and a lower bound of the bisectional curvature, we obtain
\begin{equation}
\frac{\Delta'\vert \nabla \varphi \vert ^{2}}{\vert \nabla \varphi \vert ^{2}}\geq \frac{A}{\vert \nabla \varphi\vert ^{2}}-2C_4\left(1+g'^{\alpha \overline\beta }g_{\alpha \overline\beta }\right)-C_3ke^{k(C_0-\varphi )}.
\end{equation}
Reporting (15) and (17) into (13) and taking account of $\Delta '\varphi =n-g'^{\alpha \overline{\beta}}g_{\alpha \overline{\beta}}$, we see that
$$\left[ke^{k(C_0-\varphi )}-2C_4\right]g'^{\alpha \overline{\beta}}g_{\alpha \overline{\beta}}+k^{2}e^{k(C_0-\varphi )}\vert \nabla \varphi \vert '^{2}\leq C_5ke^{k(C_0-\varphi )}.$$ 
Taking $k=2C_4+1$, since $\displaystyle e^{k(C_0-\varphi
)}\geq 1$, we get
\begin{equation}
g'^{\alpha \overline{\beta}}g_{\alpha \overline{\beta}}\leq C_5 \quad \mbox{and}\quad \vert \nabla \varphi \vert '^{2}(P)\leq C_5.
\end{equation}
In a $g$-orthonormal frame such that the matrix $(\nabla _{\lambda \overline{\mu}}\varphi)$ is diagonal. We have $\displaystyle g'^{\alpha \overline{\beta}}g_{\alpha \overline{\beta}}=\sum _{\lambda =1}^n\frac{1}{1+\nabla _{\lambda\overline{\lambda}}\varphi}$ and, for each $\lambda \in \{1,\dots ,n\}$, we can write
$$1+\nabla _{\lambda\overline{\lambda}}\varphi=F\prod _{\mu \neq \lambda}\frac{1}{1+\nabla _{\mu\overline{\mu}}\varphi}\leq F\left(\frac{1}{n-1}\sum _{\lambda =1}^n\frac{1}{1+\nabla _{\lambda\overline{\lambda}}\varphi}\right)^{n-1}.$$ 
Summing over all $\lambda$ and taking account of $(2)$ and $(18)$, we get
$$n+\bigtriangleup \varphi \leq C_6|\nabla \varphi|^r$$
Finally, the inequality $\displaystyle \vert \nabla \varphi
\vert ^{2}\leq (n+\bigtriangleup \varphi ){\vert \nabla \varphi \vert'}^2$ combined with the second inequality in (18) implies that 
$$\vert \nabla \varphi \vert ^2\leq C_5C_6 \vert \nabla \varphi \vert ^r.$$
The gradient of $\varphi $ is estimated at $P$ since $r<2$. So, the proof of theorem 1 is complete.

\vskip8mm

\section{Proof of theorem 2}

\vskip4mm

\noindent Using our gradient estimate, we can adapt the calculations of P. Cherrier [3] to estimate the Laplacian which guarantees the equivalence between the metrics $g$ and $g'$. Namely, there exists a strictly positive constant such that
\begin{equation}
a^{-1}g\leq g'\leq a g
\end{equation}
The constant $a$ depends only on the geometry of $M$, a uniform estimates of the solution and its gradient, a bound of the $C^{2}$-norm of $F$ on the compact ${\mathscr K}$ and $\displaystyle \min _{\mathscr K}F$.

\vskip4mm

\begin{lem} Let $\varphi \in C^{3}(M)$ be an admissible
solution of $(3)$ such that $\|\varphi\|_{\infty}\leq C_0$ and $|\nabla \varphi|\leq C_0$. Denote by $f=\log F$,
$$\psi ^{2}=g'^{\alpha \overline{\beta 
}}g^{a\overline{b}}\nabla _{\alpha  a}\varphi \nabla  _{\overline{\beta b}}\varphi \quad \mbox{and}\quad \theta ^{2}=g'^{\alpha \overline{\beta  }}g'^{a\overline{b}}g^{\lambda \overline{\mu }}\nabla _{\alpha  \overline{b}\lambda }\varphi \nabla _{\overline{\beta }a\overline{\mu }}\varphi .$$
Then there exists a constant $c$ such that

\begin{equation}
\begin{array}{rcc} \Delta '\psi ^{2}&\geq&2\mathrm{Re}\left(g'^{\alpha \overline{\beta  }}g^{a\overline{b}}\nabla _{\alpha  a}\varphi \nabla _{\overline{\beta  b}}f\right)-g'^{\alpha \overline{\sigma  }}g'^{\rho \overline{\beta  }}g^{a\overline{b}}\nabla _{\alpha a}\varphi \nabla  _{\overline{\beta  b}}\varphi \nabla  _{\rho \overline{\sigma }}f\\[3mm]&&-c(1+\theta )(1+\psi ^{2})+\Vert  E\Vert '^{2}+K_1,\end{array}
\end{equation}
where, in a $g'$-orthonormal frame, $\displaystyle K_1=\sum_{\lambda ,\alpha ,a}g^{a\overline{a}}\vert \nabla  _{\overline{\lambda}\alpha a}\varphi \vert ^{2}$ and the components of the tensor E are given by

\begin{equation}
E_{\lambda  \alpha  a}=\nabla  _{\lambda \alpha  a}\varphi -\sum _{\rho }\nabla  _{\lambda {\overline{\rho }}\alpha }\varphi \nabla _{\rho a}\varphi .
\end{equation}
The constant $c$ depends on $C_0$, $C_1$, the constant $a$ in $(19)$ as well as a uniform estimate of the norms of the curvature tensor ${\cal R}$, the torsion tensor ${\cal R}$ and the tensor $\nabla  {\cal T}$.
\end{lem}

\vskip6mm

\noindent\textit{Proof.} Let us compute $\displaystyle \Delta '\psi ^{2}$,
\begin{equation}
\Delta '\psi ^{2}=\sum _{i=1
}^6K_i,
\end{equation}
where

\begin{enumerate}

\item[] $K_1=g'^{\lambda \overline{\mu }}g'^{\alpha \overline{\beta}}g^{a\overline{b}}\nabla _{\overline{\mu }\alpha a}\varphi \nabla _{\lambda \overline{\beta b}}\varphi \; $ and $\; K_2=g'^{\lambda \overline{\mu }}g'^{\alpha \overline{\beta}}g^{a\overline{b}}\nabla _{\lambda \alpha a}\varphi \nabla_{\overline{\mu \beta b}}\varphi $

\vskip2mm

\item[] $\displaystyle K_3=g'^{\lambda  \overline{\mu }}g^{a\overline{b}}\left(g'^{\alpha\overline{\delta}}g'^{\gamma \overline{\sigma }}g'^{\rho  \overline{\beta }}+g'^{\alpha \overline{\sigma }}g'^{\rho 
\overline{\delta }}g'^{\gamma\overline{\beta }}\right)\nabla _{\lambda \overline{\delta }\gamma }\varphi\nabla _{\overline{\mu }\rho \overline{\sigma }}\varphi \nabla _{\alpha a}\varphi\nabla _{\overline{\beta b}}\varphi $

\item[] $\begin{array}{r}\displaystyle K_4=-g'^{\lambda  \overline{\mu }}g'^{\alpha\overline{\sigma}}g'^{\rho\overline{\beta }}\Big[\Big(\nabla _{\lambda \overline{\sigma}\rho}\varphi \nabla _{\overline{\mu }\alpha a}\varphi +\nabla _{\overline{\mu }\rho \overline{\sigma }}\varphi \nabla _{\lambda\alpha a}\varphi \Big)\nabla _{\overline{\beta b}}\varphi\\ [4mm]\displaystyle +\left(\nabla  _{\lambda \overline{\sigma  }\rho  }\varphi \nabla_{\overline{\mu \beta b}}\varphi +\nabla _{\overline{\mu }\rho \overline{\sigma }}\varphi \nabla _{\lambda \overline{\beta  b}}\varphi \right)\nabla _{\alpha a}\varphi \Big]\end{array}$

\item[] $\displaystyle K_5=g'^{\lambda \overline{\mu }}g'^{\alpha \overline{\beta}}g^{a\overline{b}}\left(\nabla _{\lambda \overline{\mu }\alpha  a}\varphi\nabla _{\overline{\beta  b}}\varphi +\nabla _{\alpha a}\varphi \nabla_{\lambda \overline{\mu \beta b}}\varphi \right)$

\item[] $\displaystyle K_6=-g'^{\lambda \overline{\mu  }}g'^{\alpha \overline{\sigma  }}g'^{\rho \overline{\beta }}g^{a\overline{b}}\nabla _{\lambda  \overline{\mu  }\rho  \overline{\sigma}}\varphi \nabla  _{\alpha  a}\varphi \nabla  _{\overline{\beta  b}}\varphi$
\end{enumerate}

\vskip2mm

\noindent We will note that $\displaystyle U(\varphi) \simeq V(\varphi)$, and say that $U$ and $V$ are equivalent, if there exists a positive constant $c$, depending only on the quantities of the lemma, such that 
$$\vert U(\varphi )-V(\varphi )\vert \leq c(1+\theta 
)(1+\psi ^{2}).$$

\vskip2mm

\noindent For a sake of simplicity, we will often express the tensors in $g'$-orthonormal frames and make use of the convention of summing over repeated indexes.

\vskip3mm

\noindent (i) Study of $K_3$. In a $g'$-orthonormal frame, if we set
$$A=g^{a\overline{a}}\nabla  _{\lambda  \overline{\alpha }\gamma  }\varphi \nabla  _{\overline{\lambda  }\rho  \overline{\gamma }}\varphi \nabla 
_{\alpha  a}\varphi \nabla  _{\overline{\rho  a}}\varphi ,$$
and 
$$B=g^{a\overline{a}}\nabla  _{\lambda  \overline{\rho
}\gamma  }\varphi \nabla  _{\overline{\lambda  }\rho  \overline{\alpha }}\varphi \nabla  _{\alpha  a}\varphi \nabla  _{\overline{\gamma  a}}\varphi
$$
then we see that
\begin{equation}
K_3=A+B.
\end{equation}

\vskip2mm

\noindent (ii) Study of $K_4$. In a $g'$-orthonormal frame, setting
$$C=g^{a\overline{a}}\nabla  _{\lambda  \overline{\alpha
 }a}\varphi \nabla _{\overline{\lambda  }\alpha  a}\varphi \nabla  _{\overline{\rho a}}\varphi $$
and
$$D=g^{a\overline{a}}\nabla  _{\lambda \overline{\alpha
 }\rho }\varphi \nabla  _{\overline{\lambda  \rho a}}\varphi \nabla _{\alpha  a}\varphi ,$$
we see that
\begin{equation}
K_4=-(C+{\overline{C }})+(D+{\overline{D }}).
\end{equation}

\vskip2mm

\noindent (iii) Study of $K_5$. In view of the fact that $\vert \nabla  _{\lambda  \overline{\mu  }\alpha  a}\varphi -\nabla _{\overline{\mu }\lambda  \alpha  a}\varphi \vert  \leq c\psi $, we can write
\begin{equation}
K_5\simeq K'_5+\overline{K'_5}
\end{equation}
where
$$K'_5=g'^{\lambda  \overline{\mu  }}g'^{\alpha 
\overline{\beta}}g^{a\overline{b}}\nabla _{\lambda  \overline{\mu }\alpha  a}\varphi\nabla _{\overline{\beta b}}\varphi .$$
Differentiating twice equation (1) satisfied by $\varphi
$,
$$g'^{\lambda  \overline{\mu }}\nabla 
_{a\lambda \overline{\mu }}\varphi =\nabla  _{a}f$$
$$g'^{\lambda  \overline{\mu }}\nabla _{\alpha  a\lambda \overline{\mu  }}\varphi =\nabla  _{\alpha 
a}f+g'^{\lambda  \overline{\rho }}g'^{\nu  \overline{\mu  }}\nabla _{\alpha  \nu  \overline{\rho  }}\varphi
\nabla  _{a\lambda  \overline{\mu }}\varphi .$$
Note that
$$g'^{\lambda  \overline{\mu  }}\nabla  _{\lambda
\overline{\mu  }\alpha  a}\varphi =g'^{\lambda \overline{\mu 
}}(\nabla _{\lambda  \overline{\mu  }\alpha a}\varphi -\nabla  _{\alpha a\lambda \overline{\mu  }}\varphi )+\nabla _{\alpha  a}f+g'^{\lambda \overline{\rho }}g'^{\nu \overline{\mu }}\nabla  _{\alpha  \nu \overline{\rho }}\varphi\nabla _{\alpha \lambda  \overline{\mu }}\varphi .$$
So, by straightforward calculation using only the definition of the covariant derivative, we see that
$$\vert  \nabla  _{\alpha  a\lambda  \overline{\mu
 }}\varphi -\nabla _{\lambda  \overline{\mu  }\alpha  a}\varphi \vert \ \leq \ C(1+\theta 
)(1+\psi ) .$$
Thus
$$K'_5\simeq g'^{\alpha  \overline{\beta }}g^{a\overline{b}}\Big\{\nabla  _{\alpha  a}f+g'^{\lambda  \overline{\rho 
}}g'^{\nu  \overline{\mu}}\nabla  _{a\lambda  \overline{\mu  }}\varphi \nabla  _{\alpha  \nu \overline{\rho  }}\varphi \Big \}\nabla _{\overline{\beta  b}}\varphi .$$ 
In this equivalence, according to relations (5) and (6), the term of type $(\nabla ^{3}\varphi )^{2}\nabla  ^{2}\varphi $ is equivalent to $C$. So taking account of (28), we have
\begin{equation}
K_5\simeq C+{\overline{C }}+\mathrm{Re}\left(g^{a{\overline{a }} }\nabla _{{\overline{\alpha  a }} }\varphi \nabla  _{\alpha  a}f\right)
\end{equation}

\vskip2mm

\noindent (iv) Study of $K_6$. We differentiate twice the equations $(1)$ to eliminate the forth order derivatives. We obtain

\begin{equation}
\begin{array}{c} K_6+g'^{\alpha 
\overline{\sigma  }}g'^{\rho \overline{\beta  }}g^{a\overline{b}}(\nabla _{\rho  \overline{\sigma }}f+g'^{\lambda  \overline{\delta  }}g'^{\gamma 
\overline{\mu  }}\nabla _{\rho  \overline{\delta  }\gamma  }\varphi \nabla _{\overline{\sigma }\lambda  \overline{\mu  }}\varphi )\nabla  _{\alpha a}\varphi \nabla _{\overline{\beta  b}}\varphi \\[3mm]\displaystyle =g'^{\lambda  \overline{\mu  }}g'^{\alpha \overline{\sigma  }}g'^{\rho \overline{\beta  }}g^{a\overline{b}}(\nabla _{\rho  \overline{\sigma }\lambda  \overline{\mu  }}\varphi -\nabla _{\lambda  \overline{\mu  }\rho  \overline{\sigma  }}\varphi )\nabla _{\alpha  a}\varphi \nabla _{\overline{\beta  b}}\varphi .\end{array}
\end{equation}

\noindent Relying on the definition of the covariant derivative, standard computations give the following identity
\begin{equation}
\begin{array}{ccl} \nabla _{\alpha \overline{\beta \mu }\lambda }\varphi -\nabla _{\overline{\mu }\lambda \alpha\overline{\beta }}\varphi &=&T^{\gamma }_{\lambda \alpha }\nabla_{\overline{\mu }\gamma  \overline{\beta }}\varphi+\Big(R^{\gamma}_{\lambda \overline{\mu }\alpha }+\nabla _{\overline{\mu }}T^{\gamma}_{\lambda  \alpha }\Big)\nabla _{\gamma  \overline{\beta }}\varphi\\[4mm] && \displaystyle +T^{\overline{\gamma}}_{\overline{\mu}\overline{\beta}}\nabla _{\lambda \overline{\gamma }\alpha}\varphi +\Big(R^{\overline{\gamma }}_{\overline{\beta \mu }\alpha }+\nabla_{\alpha }T^{\overline{\gamma}}_{\overline{\mu \beta }}\Big)\nabla _{\lambda \overline{\gamma }}\varphi.\; \end{array}
\end{equation}
It follows from this relation that the right hand term in (24) in equivalent to zero. On the other hand, in a $g'$-orthonormal frame, the term of type $(\nabla ^{2}\varphi )^{2}(\nabla ^{3}\varphi )^{2}$ in the left side of (24) can be written as follows
$$g^{a\overline{a}}\nabla _{\gamma \overline{\rho
 }\lambda  }\varphi \nabla  _{\overline{\alpha  }\rho  \overline{\lambda }}\varphi \nabla  _{\alpha  a}\varphi \nabla  _{\overline{\gamma a}}\varphi .$$
And, taking account of (5) and (6), this term differs from $B$ by $\displaystyle C(1+\theta )\psi ^{2}$. Therefore, 
\begin{equation}
\vert K_6+g'^{\alpha  \overline{\sigma  }}g'^{\rho 
\overline{\beta }}g^{a\overline{b }}\nabla _{\alpha  a
}\varphi \nabla _{\overline{\beta b } }\varphi \nabla  _{\rho \overline{\sigma }}f+B\vert \leq c(1+\theta )\psi ^{2}.
\end{equation}

\vskip0mm

\noindent (v) Finally, reporting (23), (26), (27) and (29) into (22), we get
$$\Delta '\psi ^{2}\simeq K_1+K_2+A+2\mathrm{Re}\left(g^{a\overline{a}}\nabla  _{\overline{\alpha  a}}\varphi
\nabla  _{\alpha a}f-D\right)-g^{a\overline{a}}\nabla  _{\alpha a}\varphi \nabla _{\overline{\rho  a}}\varphi \nabla  _{\rho \overline{\alpha  }}f.$$
That is (20) since the definition (21) implies that
$$\Vert E\Vert '^{2}=g'^{\lambda \overline{\mu }}g'^{\alpha \overline{\beta }}g^{a\overline{b}}E_{\lambda \alpha a}E_{\overline{\mu \beta b}}=K_2-(D+\overline{D})+A.$$

\vskip2mm

\noindent\textbf{Proof of theorem 2.} Let $k$ and $\ell$ be two strictly positives real. Set $\displaystyle J=k\vert  \nabla  \varphi \vert ^{2}+\ell(n+\Delta \varphi )$ and consider
$$\Gamma (\varphi )=\psi exp(e^{J}) .$$
Since $M$ is compact, $\Gamma (\varphi)$ attains its maximum at a point $P\in  M$. Suppose that $\psi (P)\geq 1$, otherwise we are done. Writing that the Laplacien of $\displaystyle \log\Gamma (\varphi )=\frac{\log \psi ^2}{2}+e^J$ in the metric $g'$ is non negative gives
\begin{equation}
\frac{\Delta '\psi ^2}{2\psi ^2}-\frac{\vert \nabla \psi ^2\vert '^2}{2\psi ^4}+e^{J}\vert \nabla J\vert '^2+e^{J}\Delta 'J\leq 0.
\end{equation}

\vskip2mm

\noindent We have to bound from below each term in (30). The expansion of the term involving $f$ in (20) gives a positive constant $C_3$ such that
\begin{equation}
\Delta '\psi ^2\geq -C_3\left(1+\theta ^{2}+\psi^{4}\right)+(K_1+\Vert  E\Vert'^{2})+2\mathrm{Re}\left( \frac{\partial  f_p}{\partial p^{\alpha  }}\nabla  _{\alpha  }\psi ^{2}\right).
\end{equation}

\vskip2mm

\noindent Now let us expand $\displaystyle \nabla \psi ^{2}$ in a $g$-orthonormal frame such that $g'$ is diagonal. In view of (21), we can write
$$\nabla _{\lambda }\psi ^2=g^{a\overline{a}}E_{\lambda\alpha  a}\nabla _{\overline{\alpha  a}}\varphi +g^{a\overline{a}}\nabla _{\alpha a}\varphi \nabla _{\lambda  \overline{\alpha a}}\varphi := A_{\lambda }+B_{\lambda },$$

\noindent where $\displaystyle A_{\lambda }=g^{a\overline{a}}E_{\lambda  \alpha a}\nabla  _{\overline{\alpha  a}}\varphi $ and $\displaystyle B_{\lambda }=g^{a\overline{a}}\nabla  _{\alpha  a}\varphi \nabla _{\lambda \overline{\alpha  a}}\varphi $. Next, with
$$\vert A\vert  ^{2\ }=\ \sum _{\lambda  }\vert  A_{\lambda }\vert  ^{2}\quad \mbox{et}\quad\vert  B\vert  ^{2}=\sum _{\lambda  }\vert B_{\lambda  }\vert  ^{ 2} ,$$
for all $\varepsilon  >0$, we have
$$\displaystyle \vert \nabla \psi ^{2}\vert '^{2}=\sum _{\lambda }\vert A_{\lambda }+B_{\lambda }\vert 
^{2}\leq (1+\varepsilon )\vert A\vert ^{2}+(1+\varepsilon 
^{-1})\vert B\vert ^{2}.$$
However, the inequality of Cauchy implies that
$$\vert A\vert ^2\leq \ \psi ^{2}\Vert E\Vert '^2 \quad\mbox{and}\quad \vert B\vert  ^2\leq  \psi ^{2}K_1.$$
So

\begin{equation}
\frac{\vert \nabla \psi ^{2}\vert '^2}{\psi ^2}\;
\leq \; (1+\varepsilon )\Vert E\Vert '^{2}+(1+\varepsilon  ^{-1})K_1.
\end{equation}
On the other hand, since $\log\Gamma (\varphi )$ is stationary at $P$, we have
\begin{equation}
\frac{\nabla \psi ^2}{2\psi ^2}=-e^{J}\nabla J.
\end{equation}
Set $\displaystyle \delta  =2^{-1}e^{-J}$. Relation (33) implies that
$$e^{J}\vert \nabla  J\vert '^2=e^{-J}(e^{2J}\vert \nabla J\vert '^2)=e^{-J}\frac{\vert \nabla \psi ^2\vert '^2}{4\psi ^4}=\delta \frac{\vert \nabla \psi ^2\vert '^2}{2\psi ^4}.$$
So, we can write
$$\frac{\vert \nabla \psi ^2\vert '^2}{2\psi ^4}-e^{J}\vert \nabla J\vert '^{2}=(1-\delta )\frac{\vert \nabla \psi
^{2}\vert '^2}{2\psi ^4}$$
and, in view of (32),

$$
\frac{\vert \nabla \psi ^{2}\vert '^2}{2\psi ^4}-e^{J}\vert  \nabla J\vert '^2\leq (1-\delta)\frac{(1+\varepsilon )\Vert E\Vert '^{2}+(1+\varepsilon  ^{-1})K_1}{2\psi ^2}.
$$

\noindent Choose $\displaystyle \varepsilon  =\frac{\delta}{1-\delta}$ in (32). Since $\displaystyle (1-\delta )(1+\epsilon )=1$ and $\displaystyle (1-\delta)(1+\epsilon ^{-1})=\frac{1-\delta}{\delta}=2e^{J}-1$, we obtain
\begin{equation}
\frac{\vert \nabla \psi ^{2}\vert '^2}{2\psi ^4}-e^{J}\vert  \nabla J\vert '^2\leq \frac{\Vert E\Vert 
'^{2}+(2e^{J}-1)K_1}{2\psi
^2}.
\end{equation}

\noindent We pursue by the study of the last term in (30). The relation in Lemma 1 gives,
$$\Delta '\vert \nabla  \varphi \vert ^{2}\geq \ -C_4+\psi ^{2}+2\mathrm{Re}\left(\frac{\partial f_p}{\partial p^{\alpha }}\nabla _{\alpha }\vert \nabla \varphi \vert ^{2}\right)$$

\noindent Where the constant $C_4$ is under control. Now, we differentiate twice equation (1),

$$\nabla _{\lambda  }f=g'^{\alpha  \overline{\beta
 }}\nabla  _{\lambda  \alpha \overline{\beta  }}\varphi $$ 
and
$$\nabla _{\lambda \overline{\mu}}f=g'^{\alpha\overline{\beta  }}\nabla _{\lambda \overline{\mu }\alpha\overline{\beta  }}\varphi -g'^{\alpha \overline{\sigma }}g'^{\rho\overline{\beta }}\nabla _{\overline{\mu }\rho \overline{\sigma }}\varphi \nabla _{\lambda \alpha \overline{\beta }}\varphi .$$
So
$$\begin{array}{ccl}\Delta f&=&g^{\lambda \overline{\mu }}g'^{\alpha \overline{\beta }}\nabla _{\lambda \overline{\mu }\alpha\overline{\beta }}\varphi -g'^{\alpha \overline{\sigma  }}g'^{\rho \overline{\beta }}g^{\lambda \overline{\mu }}\nabla  _{\overline{\mu }\rho \overline{\sigma }}\varphi
\nabla _{\lambda \alpha \overline{\beta }}\varphi \\[3mm] &=&g'^{\alpha\overline{\beta}}(n+\Delta \varphi)+g^{\lambda  \overline{\mu  }}g'^{\alpha \overline{\beta }}\left(\nabla _{\lambda \overline{\mu }\alpha\overline{\beta }}\varphi-\nabla _{\alpha\overline{\beta }\lambda \overline{\mu }}\varphi\right)-g'^{\alpha \overline{\sigma }}g'^{\rho
\overline{\beta }}g^{\lambda \overline{\mu }}\nabla  _{\overline{\mu }\rho \overline{\sigma }}\varphi
\nabla _{\lambda \alpha \overline{\beta}}\varphi .\end{array}$$
Computing $\Delta f$ and using the identity (28), we see that
$$\Delta '(n+\Delta\varphi )\geq -C_5\psi ^2+\theta  ^{2}+2\mathrm{Re}\left(\frac{\partial f_p}{\partial  p^{\alpha}}\nabla _{\alpha }(n+\Delta \varphi )\right),$$
where $C_5$ is a constant under control. Therefore

\begin{equation}
\Delta 'J\; \geq \; -C_4k+\ell\theta ^2+(k-C_5\ell)\psi ^{2}+2\mathrm{Re}\left( \frac{\partial f_p}{\partial  p^{\alpha } }\nabla _{\alpha }J\right).
\end{equation}
Finally, injecting (31), (34) and (35) into (30) and making use of (33), we arrive at the following inequality
\begin{equation}
-\frac{C_3(1+\theta ^2+\psi ^4)}{2\psi ^2}-C_4e^J+\ell e^J\theta ^2+(k-C_5\ell)e^J\psi ^2-\frac{e^{J}-1}{\psi
^2}K_1\leq 0.
\end{equation}
By straightforward computations, we show that $\displaystyle \nabla _{\alpha \overline{\lambda}a}\varphi -\nabla _{\overline{\lambda }\alpha  a}\varphi =R^{b}_{a\alpha\lambda }\nabla _{b}\varphi $. So $\displaystyle K_1\leq 2\theta  ^{2}+C_6$, where $C_6$ is a controlled constant.
Now, since $\psi (P)\geq 1$, inequality (36) takes the following form
\begin{equation}
(\ell e^J-C_3-2e^J)\theta ^2 +(k-C_5\ell -C_3e^{-J})e^{J}\psi
^2\leq C_3+C_4e^J+C_6.
\end{equation}
We choose $\displaystyle \ell=2+C_3$ and then $\displaystyle k=1+C_5\ell +C_3$. Since $\displaystyle e^{J}\geq 1$,
coefficients of $\theta  ^{2}$ and $\psi ^{2}$ in (37)
are respectively $\geq 0$ and $\displaystyle \geq 1$. Thus
$$\displaystyle \psi ^{2}(P)\ \leq \max\left(1,C_3+C_6+C_4e^{\|J\|_{\infty}}\right).$$
This end the proof.

  \end{document}